# A THEORY OF STOCHASTIC INTEGRATION FOR BOND MARKETS

By M. De Donno and M. Pratelli

*Università di Pisa*

We introduce a theory of stochastic integration with respect to a family of semimartingales depending on a continuous parameter, as a mathematical background to the theory of bond markets. We apply our results to the problem of super-replication and utility maximization from terminal wealth in a bond market. Finally, we compare our approach to those already existing in literature.

**1. Introduction.** In the models for the bond market it is usually assumed that a continuum of basic securities is available to the investor. This gives rise to the problem of what exactly should be meant by the word "portfolio" in this setting. For this reason, Björk, Di Masi, Kabanov and Runggaldier [4] modeled the zero coupon bonds price process as a stochastic process with values in the space $C([0, \infty))$ of continuous functions on the time interval $[0, \infty)$. Then, they constructed a stochastic integral with respect to such a process, where the integrand process (i.e., the mathematical representation of a self-financing strategy) takes values in the dual of $C([0, \infty))$, that is, the set of Radon measures on $[0, \infty)$. In this way, they define a notion of "infinite-dimensional" portfolio as a portfolio based on a measure-valued strategy and a contingent claim is called *attainable* when its value coincides with the final value of a measure-valued portfolio. This approach has, however, some drawbacks. For instance, if one wishes to extend the results on completeness from the stock market, one discovers that the *uniqueness* of the martingale probability is not equivalent to *completeness* in the usual sense, but to a weaker condition, called *approximate completeness*: every (sufficiently integrable) contingent claim can be approximated by a sequence of attainable claims, but it may not be attainable. We suggested a different approach in [10], by making use of a theory on stochastic integration with respect to cylindrical locally square integrable martingales, developed by Mikulevicius and Rozovskii [20, 21]: there, we showed that measure-valued processes are not sufficient to describe all financial portfolios. Nevertheless, this type of approach is limited to the martingale case, which means that it is necessary to work under an equivalent martingale measure, whereas there are some questions, such as super-replication or utility maximization, which need to be posed under the original measure.

The aim of this paper is to introduce a theory of stochastic integration with respect to a family of semimartingales indexed by a continuous parameter $x \in I$ (where $I$ is a locally compact subset of $\mathbb{R}$) and to characterize the class of integrands. We start from defining as *simple integrands* the linear combination of Dirac measures, so that the integral depends only on a finite number of semimartingales. Going to the limit (in a properly chosen topology), we define the class of *generalized integrands*, which turns out to be predictable processes with values on the set of not necessarily bounded functionals on $\mathbb{R}^I$.

This approach, which, in our opinion, is the most natural from a mathematical point of view, makes sense also from the financial point of view. Indeed, the simple integrands are the mathematical representation of *real* portfolios, which are based on a finite number of bonds: it is then natural to take as generalized portfolios the limit of *real* portfolios.

We introduce our definitions in Section 2, where we show that the integral is well defined. Moreover, we provide an infinite-dimensional analog of a result due to Mémin [18]: the limit of stochastic integrals (for the semimartingale topology) is still a stochastic integral. This section is the extension to the continuous case of the results obtained in [11], concerning the stochastic integration with respect to a sequence of semimartingales.

In Section 3 we study the problems of *super-replication* and *utility maximization* in the bond market. We adapt the results of Delbaen and Schachermayer [12] and Kramkov and Schachermayer [17] by making use of the theory developed in Section 2. We remark that these results cannot be obtained if we consider only measure-valued strategies, as it was observed by Pham [22].

In Section 4 we compare our theory with other approaches existing in literature, which are also based on infinite-dimensional stochastic integration. At first we concentrate on the Hilbert space approach, according to which the family of bond price processes is seen as a process with values in a Hilbert space of functions, instead of a family of semimartingales dependent on a parameter. In particular, we analyze the recent paper by Carmona and Tehranchi [6].

Then, we discuss connections with the paper by Björk, Di Masi, Kabanov and Runggaldier [4]: in particular, we show that the stochastic integral constructed in that paper can be seen as a particular case of our definition of generalized integral.

We point out that *generalized strategies* in the bond market models are an idealization, impossible to achieve in the real word: they can only be approximated by *real* portfolios. However, in some circumstances, the language of generalized integrands allows a simpler analysis of the problem and a better understanding of the real bond market: for instance, in the problem of *utility maximization* we can give an explicit form of the optimal (idealized) solution by using generalized integrands, while the optimal solution might not exist if we consider only portfolios based on a finite number of assets.

**2. Definitions and general results.** Let $(\Omega, \mathcal{F}, (\mathcal{F}_t)_{0 \leq t \leq T}, \mathbb{P})$ be a filtered probability space, which satisfies the usual hypotheses. Let $\mathcal{S}(\mathbb{P})$ be the space of real-valued semimartingales. We consider $\mathcal{S}(\mathbb{P})$ equipped with the topology introduced by Emery [15]: $\mathcal{S}(\mathbb{P})$ is a complete metric space. It is worth it to recall that $\mathcal{S}(\mathbb{P})$ is invariant with respect to a change to an equivalent probability measure. We denote by $\| \cdot \|_{\mathcal{S}(\mathbb{P})}$ the quasinorm introduced by Mémin ([18], Section II), which induces the Emery topology on the set of semimartingales. Moreover, we denote by $\mathcal{M}^2(\mathbb{P})$ the Banach space of square integrable martingales with the norm $\|M\|_{\mathcal{M}^2} = \|M_T\|_{L^2(\mathbb{P})}$. Finally, $\mathcal{A}(\mathbb{P})$ is the Banach space of predictable processes with finite variation, whose variation is integrable, with the norm $\|B\|_{\mathcal{A}} = \|V(B)_T\|_{L^1(\mathbb{P})}$ (with the notation of Mémin [18], $V(B)$ denotes the variation of $B$). Many of our results rely on Theorem II.3 of [18]:

*Let $S^n$ be a Cauchy sequence in $\mathcal{S}(\mathbb{P})$; there exist a probability measure $\mathbb{P}^*$, equivalent to $\mathbb{P}$, such that $d\mathbb{P}^*/d\mathbb{P}$ is bounded and a subsequence (still denoted by $S^n$) such that $S^n = M^n + B^n$, where $M^n$ is a Cauchy sequence in $\mathcal{M}^2(\mathbb{P}^*)$ and $B^n$ is a Cauchy sequence in $\mathcal{A}(\mathbb{P}^*)$.*

With the notation introduced by Mémin [18], we will say that $S^n$ converges in $\mathcal{M}^2 \oplus \mathcal{A}(\mathbb{P}^*)$. Since this result holds only for processes defined on a finite horizon, we limit our study to the case of processes defined on a finite time set $[0, T]$. However, we point out that the choice of a finite time interval is not restrictive: indeed, all the results can be extended to the time interval $[0, +\infty[$ by localization.

Let $\mathbf{S} = (S^x)_{x \in I}$ be a family of semimartingales, where $I$ is a locally compact subset of $\mathbb{R}$; in fact, it can be replaced with any locally compact separable metric space. The infinite-dimensional semimartingale $\mathbf{S}$ can be seen as a mapping from a metric space to the space of semimartingales, endowed with the Emery topology. We make the following assumption:

ASSUMPTION 2.1. *The mapping $\mathbf{S} : I \to \mathcal{S}(\mathbb{P})$ defined by $\mathbf{S}(x) = S^x$ is continuous.*

Our goal is to give a good definition of a stochastic integral with respect to the infinite-dimensional semimartingale $\mathbf{S}$. Most of the definitions and results in this section will be an extension, to the "continuous" case, of the definitions and results given in [11] for the case of a sequence of semimartingales.

In [11], we called negligible [with respect to a sequence of semimartingales $(X^i)_{i\geq 1}$] a predictable set $C \subset \Omega \times [0, T]$ such that $\int h \, dX^i = h \cdot X^i = 0$, for every $i$ and for every bounded predictable process $h$, which is zero on the complement of $C$. We also proved that there exists an increasing predictable process $A$, such that a set is negligible (according to the above definition) if and only if it is negligible with respect to the measure $d\mathbb{P} \, dA$. In the present setting we introduce the following definition:

DEFINITION 2.1.   Let $\{x_n\}_{n\geq 1}$ be a countable dense subset in $I$: we will call *negligible* [with respect to the family $(S^x)_{x\in I}$] a predictable set $C$ which is negligible with respect to the sequence $(S^{x_n})_{n\geq 1}$.

The next lemma shows that this notion is well defined, in the sense that it does not depend on the choice of the dense subset.

LEMMA 2.1.   *For a predictable set $C$, the following conditions are equivalent:*

(i) *$C$ is negligible with respect to the sequence of semimartingales $(S^{x_n})_{n\geq 1}$ (where $\{x_n\}_{n\geq 1}$ is a countable dense subset in $I$).*

(ii) *For every $x \in I$, for every bounded predictable process $h$ which is zero on the complement of $C$, we have that $(h \cdot S^x) = 0$.*

PROOF.   It is trivial that (ii) implies (i). Assume then that (i) holds. Condition (ii) is obvious if $x$ belongs to the dense set considered in (i); otherwise, it follows from Assumption 2.1. In details, fix $x \notin \{x_n\}_{n\geq 1}$ and let $h$ be a predictable process such that $|h(\omega, t)| \leq M$ for all $(\omega, t)$ and $h(\omega, t) = 0$ for $(\omega, t) \notin C$.

Then, recalling that $h \cdot S^{x_n} = 0$ for all $n$ because of (i), we have that

$$\|h \cdot S^x\|_{\mathcal{S}(\mathbb{P})} = \|h \cdot S^x - h \cdot S^{x_n}\|_{\mathcal{S}(\mathbb{P})} \leq M\|S^x - S^{x_n}\|_{\mathcal{S}(\mathbb{P})}.$$

Moreover, by Assumption 2.1, for every $\varepsilon > 0$, there exists $n_\varepsilon$ such that $\|S^x - S^{x_{n_\varepsilon}}\|_{\mathcal{S}(\mathbb{P})} < \varepsilon$. Since $\varepsilon$ can be chosen arbitrarily small, the claim follows.   □

From now on, when we say that a property holds a.s., we mean that it holds on a set whose complementary is negligible with respect to **S**.

In order to define a stochastic integral with respect to **S**, we begin by introducing the definition of simple integrand.

DEFINITION 2.2. A *simple integrand* (*simple strategy*) is a process $H$ of the form

$$(2.1) \qquad H_{\omega,t} = \sum_{i \leq n} h^i_{\omega,t} \delta_{x_i},$$

where $h^i$ are predictable bounded processes, whereas the $\delta_{x_i}$ denote the Dirac deltas at points $x_i \in I$.

The stochastic integral of a simple integrand with respect to $\mathbf{S}$ is naturally defined by the formula

$$(H \cdot \mathbf{S})_t = \int_0^t H_s \, d\mathbf{S}_s = \int_0^t \sum_{i \leq n} h^i_s \, dS^{x_i}_s.$$

We observe that $\mathbf{S}$ is a process with values in $E = \mathbb{R}^I$. We consider on $E$ the product topology. Hence, its dual space $E'$ is the space of finite combination of Dirac measures on $I$. It is easy to recognize that a simple integrand $H$ is a process which takes values in $E'$.

We note that our class of simple integrands coincides with the class of simple integrands introduced by Mikulevicius and Rozovskii ([20], page 142) for the case of locally square integrable cylindrical martingales in a quasicomplete and locally convex space $E$.

REMARK 2.1. If one has in mind a theory of stochastic integration for the bond market, it might seem a more convenient choice to take as value space of the infinite-dimensional semimartingale the set $\mathcal{C}(I)$ of the continuous functions from $I$ to $\mathbb{R}$. This is, in fact, a common approach in the existing literature (see, e.g., [4] and Section 4 for a comparison with their results). In this case, taking as value space of the simple integrands the dual space $(\mathcal{C}(I)')$, one should consider as simple integrands the processes with values in the set of Radon measures with compact support, but not necessarily with finite support. This seems counterintuitive to us from the financial point of view, since it should be quite evident that the simple strategies (the strategies of the real world) involve only a finite number of assets. Moreover, also from the mathematical point of view, the pointwise continuity does not imply (without any further condition) that Assumption 2.1 holds. At the same time, this assumption seems to us the right assumption to work with.

It is clear that the set of simple integrands is too small. That is, we need more complicated (though unrealistic) strategies, if we wish to address, at least theoretically, to questions as super-replication and utility maximization in a market where a continuum of securities is available.

It was already seen by Métivier [19], for the case of martingales with values in a Hilbert space $K$, that the integrands do not necessarily take values in $K'$. In particular, they take values in the set of unbounded functionals on $K$.

DEFINITION 2.3. A *not-necessarily continuous* (*unbounded*) *functional* on $E$ is a linear functional $k$ whose domain $\mathrm{Dom}(k)$ is a subspace of $E$.

We denote by $\mathcal{U}$ the set of (not-necessarily continuous) functionals on $E$.

The space $E'$ is of course contained in $\mathcal{U}$. Furthermore, we note that $\mathcal{U}$ contains also the set of Radon measures on $I$, which we denote by $\mathcal{M}(I)$: if $\mu \in \mathcal{M}(I)$, then $\mathrm{Dom}(\mu) \supset \mathcal{C}(I)$.

DEFINITION 2.4. We will say that a sequence $(k^n)$ in $E'$ *converges to* $k \in \mathcal{U}$ if $\lim_n k^n(f) = k(f)$, for all $f$ in $\mathrm{Dom}(k)$.

Note that, for a sequence $k^n$ in $E'$, it always makes sense to define the limit functional $k = \lim_n k^n$, where $\mathrm{Dom}(k) = \{f \in E : \lim_{n\to\infty} k^n(f) \text{ exists}\}$ can possibly be the trivial set $\{0\}$ (see also Remark 1 in [11]).

Since we wish to consider $\mathcal{U}$-valued processes, we introduce a notion of (weak) predictability.

DEFINITION 2.5. A process $\mathbf{H}$ with values in $\mathcal{U}$ is (*weakly*) *predictable* if for every element $f$ of $E$, the process $\mathbf{H}(f)\mathbf{1}_{\{f \in \mathrm{Dom}(\mathbf{H})\}}$ is predictable.

LEMMA 2.2. *Let* $(H^n)$ *be a sequence of simple integrands which converges to a $\mathcal{U}$-valued process* $\mathbf{H}$. *Then* $\mathbf{H}$ *is weakly predictable.*

PROOF. Let $f \in E$. The set $\{f \in \mathrm{Dom}(\mathbf{H})\} = \{(\omega, t) : \lim_n H^n_{\omega,t}(f) \text{ exists}\}$ is predictable, since each $H^n$ is predictable. Hence, the process $\mathbf{H}(f)\mathbf{1}_{\{f \in \mathrm{Dom}(\mathbf{H})\}}$ is pointwise limit of the sequence of predictable processes $H^n(f)\mathbf{1}_{\{f \in \mathrm{Dom}(\mathbf{H})\}}$, and as a consequence, is predictable as well. $\square$

We introduce as in [11], Definition 2, the notion of generalized integrand as the limit, in an appropriate sense, of simple integrands. This definition is analogous to (and, in fact, inspired by) the notion of integrable function with respect to a vector-valued measure (see, e.g., [13], Section IV.10.7).

DEFINITION 2.6. Let $\mathbf{H}$ be a $\mathcal{U}$-valued process. We say that $\mathbf{H}$ is *integrable* with respect to $\mathbf{S}$ if there exists a sequence $(H^n)$ of simple integrands such that:

(i) $H^n$ converges to $\mathbf{H}$, a.s.;
(ii) $(H^n \cdot \mathbf{S})$ converges to a semimartingale $Y$ in $\mathcal{S}(\mathbb{P})$.

We call $\mathbf{H}$ a *generalized integrand* and define $\int \mathbf{H}\,d\mathbf{S} = \mathbf{H} \cdot \mathbf{S} = Y$.
We denote by $\mathcal{L}(\mathbf{S}, \mathcal{U})$ the set of generalized integrands.

We remark that, in Definition 2.6, we do not give any condition on measurability, contrarily to what is done in the classical definition of stochastic integral. The reason is that a measurability property is implicit in (i), thanks to Lemma 2.2.

It is clear that Definition 2.6 makes sense only provided that the limit semimartingale does not depend on the approximating sequence.

PROPOSITION 2.3. *The semimartingale $Y$ of Definition 2.6 is well defined, that is, if $(H^n)$ and $(G^n)$ are sequences of simple integrands both converging to $\mathbf{H}$ and such that both $(H^n \cdot \mathbf{S})$ and $(G^n \cdot \mathbf{S})$ are Cauchy sequences in $\mathcal{S}(\mathbb{P})$, then $(H^n \cdot \mathbf{S})$ and $(G^n \cdot \mathbf{S})$ converge to the same limit.*

PROOF. We can reduce to the case of a sequence of semimartingales. Indeed, since $H^n$ and $G^n$ are simple integrands, they will have the form

$$H^n = \sum_{i \leq j(n)} h^{n,i} \delta_{x_i^n}, \qquad G^n = \sum_{i \leq k(n)} g^{n,i} \delta_{y_i^n}.$$

Therefore, if we denote by $I^*$ the set $\bigcup_{n \geq 1} (\{x_1^n, \ldots, x_{j(n)}^n\} \cup \{y_1^n, \ldots, y_{k(n)}^n\})$, we recognize that we have to deal with a stochastic integral with respect to the sequence of semimartingales $(S^x)_{x \in I^*}$. Then, by Proposition 1 in [11], the limit is unique and the stochastic integral is well defined. □

In a similar way, using Theorem 3 in [11], we can prove the following result which extends a result due to Mémin ([18], Corollary III.4) for the case of finite-dimensional semimartingales:

THEOREM 2.4. *Let $(\mathbf{H}^n)$ be a sequence of generalized integrands such that $(\mathbf{H}^n \cdot \mathbf{S})$ is a Cauchy sequence in $\mathcal{S}(\mathbb{P})$. Then, there exists a generalized integrand $\mathbf{H}$ such that $\lim_{n \to \infty} \mathbf{H}^n \cdot \mathbf{S} = \mathbf{H} \cdot \mathbf{S}$.*

The mapping $\mathcal{L}(\mathbf{S}, \mathcal{U}) \to \mathcal{S}(\mathbb{P})$ defined by $\mathbf{H} \mapsto \mathbf{H} \cdot \mathbf{S}$ is well defined, linear and invariant with respect to a change to an equivalent probability measure. However, besides these nice properties, there are some drawbacks. First of all, the integral is not linear in $\mathbf{S}$: there may exist two infinite-dimensional semimartingales $\mathbf{S}^1$ and $\mathbf{S}^2$ and a $\mathcal{U}$-valued process $\mathbf{H}$ such that $\mathbf{H}$ in $\mathcal{L}(\mathbf{S}^1, \mathcal{U}) \cap \mathcal{L}(\mathbf{S}^2, \mathcal{U})$, but $\mathbf{H} \notin \mathcal{L}(\mathbf{S}^1 + \mathbf{S}^2, \mathcal{U})$. The reason is that the integral depends on the approximating sequences and the same generalized integrand can admit two different approximating sequences according to the integrator semimartingale. Moreover, the integral is not stable for small perturbations of the semimartingale $\mathbf{S}$, as we shall show in Example 2.2. Finally, we cannot provide an extension of an important result which was proved by
Ansel and Stricker ([2], Corollary 3.5) for the case of a finite number of local
martingales. On the contrary, we can find a counterexample, that is,

*there exist a family of local martingales $\mathbf{M} = (M^x)_{x \in I}$ and a sequence of
simple integrands $(H^n)$ such that the sequence of integrals $(H^n \cdot \mathbf{M})$ con-
verges to an integral $(\mathbf{H} \cdot \mathbf{M})$ which is bounded from below, but is not a local
martingale.*

This phenomenon was already observed by De Donno and Pratelli [11],
Example 2 (which in turn is inspired to an example of Emery [15], page 496),
for the case of a sequence of semimartingales. We can modify that example
to obtain one for the "continuous" case.

EXAMPLE 2.1.   Let $M^i$ be the square integrable martingale defined as
follows:

$$M_t^i = \frac{t \wedge T_i}{i^2} - \mathbf{1}_{\{t \geq T_i\}},$$

where $(T_i)_{i \geq 1}$ is a sequence of independent random variables, such that $T_i$
is exponentially distributed with $\mathbb{E}[T_i] = i^2$. The filtration $(\mathcal{F}_t)_{t \leq T}$ is the
smallest filtration such that $T_i$ are stopping times and the usual conditions
are satisfied. Note that the sequence $M^i$ converges to 0 in $\mathcal{M}^2(\mathbb{P})$ (as $i \to \infty$),
since $\mathbb{E}[(M_T^i)^2] = \mathbb{E}[i^{-2}(T \wedge T_i)] = 1 - \exp(-T/i^2)$.

We take $I = [0, 1]$: for $x_i = 1 - i^{-1}$, we define a local martingale by $M^{x_i} =
M^i$ and by linear interpolation, we extend the mapping $x \mapsto M^x$ to the whole
$[0, 1]$ (we set $M^1 = 0$): this mapping clearly satisfies Assumption 2.1.
We take as simple integrands the processes

$$H^n = \frac{1}{n} \sum_{i \leq n} i^2 \delta_{x_i},$$

so that the integral is

$$H^n \cdot \mathbf{M} = \frac{1}{n} \sum_{i \leq n} i^2 M^i.$$

The sequence $(H^n \cdot \mathbf{M})$ is a Cauchy sequence in $\mathcal{S}(\mathbb{P})$ and converges to the
increasing process $A_t = t$. Indeed, consider the sequence of stopping times
$S_k = \inf_{m \geq k} T_m$. Using the Borel–Cantelli lemma, it can be proved that $S_k$
tends to infinity (as $k \to \infty$). In particular, the sequence $S_k \wedge T$ converges
to $T$ stationarily, namely, $S_k \equiv T$ definitely, $\mathbb{P}$-a.s. So, for fixed $\varepsilon$, there exists
some $k$ such that $\mathbb{P}(S_k \leq T) < \varepsilon$. On the stochastic interval $[[0, S_k \wedge T]]$, the
martingale $N^m = m^2 M^m$ coincides with the process $A$, for $m \geq k$. Then, if
we stop the processes $H^n \cdot \mathbf{M}$ at time $S^k$, we have that, for $n > k$,

$$(H^n \cdot \mathbf{M})^{S^k} = \frac{(N^1 + \cdots + N^k)^{S^k}}{n} + \frac{(n-k)}{n}(t \wedge S_k).$$

It is not difficult to check that the sequence $(H^n \cdot \mathbf{M})^{S_k}$ converges to $t \wedge S_k$ as $n$ tends to $\infty$: as a consequence, $(H^n \cdot \mathbf{M})$ converges to $A_t$ in $\mathcal{S}(\mathbb{P})$.

Moreover, the sequence $(H^n)$ converges (as $n \to \infty$) to the generalized integrand defined by

$$\mathbf{H}(f) = \lim_{n \to \infty} \frac{1}{n} \sum_{i \leq n} i^2 \delta_{x_i}(f) = \lim_{n \to \infty} \frac{1}{n} \sum_{i \leq n} i^2 f(x_i).$$

So, we have found a generalized integrand $\mathbf{H}$, such that the integral $(\mathbf{H} \cdot \mathbf{M})_t = t \geq 0$, but it is not a local martingale.

EXAMPLE 2.2. Let $(M^x)_{x \in [0,1]}$ and $(x_i)_{i \geq 1}$ be defined as in Example 2.1. For every $k \in \mathbb{N}$, we set $M^{k,x} = M^x$ if $x \leq x_k$ (hence, $M^{k,x_i} = M^i$ for $i \leq k$) and $M^{k,x} = 0$ if $x \geq x_{k+1}$. Then we extend $M^{k,x}$ by linear interpolation between $x_k$ and $x_{k+1}$. We observe that $\mathbf{M}^k = (M^{k,x})_{0 \leq x \leq 1}$ is a "small perturbation" of $\mathbf{M}$: indeed, for $x \leq x_k$, we have that $\|M^{k,x} - M^x\|_{\mathcal{M}^2(\mathbb{P})} = 0$, and few calculations show that, for $x > x_k$,

$$\|M^{k,x} - M^x\|_{\mathcal{M}^2(\mathbb{P})}^2 \leq 1 - e^{-T/k^2}$$

(since $\mathbb{E}[(M_T^k)^2] = 1 - \exp(-T/k^2)$). We take $H^n$ as above: for $n > k$, the integral of $H^n$ with respect to the family $\mathbf{M}^k = (M^{k,x})_{0 \leq x \leq 1}$ is given by

$$H^n \cdot \mathbf{M}^k = \frac{1}{n} \sum_{i \leq k} i^2 M^i.$$

Hence, $\mathbf{H} \cdot \mathbf{M}^k = \lim_n H^n \cdot \mathbf{M}^k = 0$ for all $k$, whereas $\mathbf{H} \cdot \mathbf{M} = A$.

## 3. Utility maximization in the bond market.

We consider a model of a bond market based on a family of semimartingales $\mathbf{P} = (P(\cdot, T))_{T \leq T^*}$; the random variable $P(t, T)$ represents the price at time $t$ of a zero coupon bond (ZCB) maturing at time $T \geq t$. We remark that $T^*$ can possibly be $\infty$: in this case, $\mathbf{P}$ is defined on the open interval $[0, \infty)$. For all the basic definitions and assumptions on models of bond markets, we mainly refer to Chapter 15 of the book of Björk [3]. In particular, we make the following (usual) assumptions to guarantee that the bond market is sufficiently rich and regular.

ASSUMPTION 3.1. (1) There exists a (frictionless) market for the ZCB for all maturities $T \leq T^*$.

(2) For each fixed $t$, the bond price $P(t, T)$ is differentiable with respect to $T$.

(3) $P(t, T) > 0$ and $P(T, T) = 1$ for all $t \leq T \leq T^*$.

We recall that the instantaneous *forward rate* at $T$, contracted at time $t$, is defined as

$$f(t,T) = -\frac{\partial P(t,T)}{\partial T};$$

the *short rate* is defined by $r(t) = f(t,t)$. The process $B_t = \exp(\int_0^t r_s\,ds)$ is the bank account and it is a strictly positive continuous process. We take this asset as a numéraire and denote by $\overline{P}(t,T)$ the discounted bond prices $P(t,T)/B(t)$.

As Björk, Di Masi, Kabanov and Runggaldier [4], we need to work with processes defined for all $t \in [0,T^*]$, while, by definition, the bond price $P(t,T)$ is given only for $t \leq T$. We then use the same trick of Björk, Di Masi, Kabanov and Runggaldier ([4], page 149), assuming that after maturity the bond is automatically transferred into the bank account. In mathematical terms, this amounts to set $P(t,T) = \exp(\int_T^t r_s\,ds) = B(t)/B(T)$ for $t > T$.

Since we want to exclude arbitrage possibilites, we make the following assumption:

Assumption 3.2. There exists an equivalent probability measure $\mathbb{P}^*$ under which the process $(\overline{P}(t,T))_{t \leq T^*}$ is a local martingale for every $T$ ($\mathbb{P}^*$ is known as an equivalent martingale measure).

We denote by $\mathcal{M}_e$ the (nonempty) set of all the equivalent martingale measures.

In order to apply the theory developed in Section 2, we need the families $\mathbf{P} = (P(t,T))_{0 \leq t, T \leq T^*}$ and $\overline{\mathbf{P}} = (\overline{P}(t,T))_{0 \leq t, T \leq T^*}$ to satisfy Assumption 2.1. We remark that this is not a consequence of Assumption 3.1. At the same time, all the models for the bond market studied in literature satisfies Assumption 2.1. Hence, it is not restrictive to assume that this condition holds.

With the notation of Section 2, we should rather write $P_t^T$ or $\overline{P}_t^T$, the latter being semimartingales indexed by a continuous parameter $T \in [0,T^*]$. Then, the theory developed in the previous section establishes exactly which are the generalized integrands (hence, the self-financing strategies) with respect to the infinite-dimensional process $\overline{\mathbf{P}} = (\overline{P}^T)_{T \leq T^*}$.

Definition 3.1. A *generalized self-financing portfolio strategy* is a pair $\pi = (V_0, \mathbf{H})$, where $V_0$ is a real number and $\mathbf{H}$ is a generalized integrand for $\overline{\mathbf{P}}$. The *discounted portfolio value process* is given by

$$\overline{V}_t^\pi = V_0 + \mathbf{H} \cdot \overline{\mathbf{P}}.$$

At this point, a short discussion is necessary on the "financial" meaning of a self-financing generalized strategy. Indeed, it is not so clear as in the finite-dimensional case either which is the relationship between a discounted and nondiscounted portfolio value or which is the investment in the bank account.

Let $\overline{V}$ be the discounted value of a generalized self-financing portfolio, generated by the generalized strategy $\mathbf{H}$. Without loss of generality, we can assume $V_0 = 0$. By definition of $\mathbf{H}$, there exists a sequence $H^n$ of simple strategies $H^n$ such that the corresponding sequence of discounted self-financing portfolios $\overline{V}_n = H^n \cdot \overline{\mathbf{P}}$ converges to $\overline{V}$ in the semimartingale topology (hence, in probability and, up to a subsequence, a.s.). If we denote by $\varphi_n$ the investment in the bank account for each of these finite-dimensional portfolios, we know that $\overline{V}_n = \varphi_n + H^n(\overline{\mathbf{P}})$ [in nondiscounted terms: $V_n = \varphi_n B + H^n(\mathbf{P})$], which entails $\varphi_n = H^n \cdot \overline{\mathbf{P}} - H^n(\overline{\mathbf{P}})$. It follows that $\lim_n \varphi_n$ exists finite (i.e., we can specify the amount invested in the bank account for the generalized strategy) if and only if $\overline{\mathbf{P}}$ belongs to the domain of $\mathbf{H}$ [or, equivalently, $\mathbf{P} \in \mathrm{Dom}(\mathbf{H})$]. In particular, we note that if $\mathbf{H}$ is a measure-valued process, the above condition is always fulfilled, since $\mathbf{P}$ takes values in $\mathcal{C}(I)$.

However, we must remember that $\mathbf{H}$ is, in any case, a theoretical strategy: it is a mathematical representation of the limit of self-financing simple strategies, the strategies of the real world, each of which has a well-determined amount invested in a bank account. Hence, an investor knows that he can buy a real portfolio which approximates the "optimal" strategy, without possibly achieving it.

Our purpose is to extend the results of De Donno, Guasoni and Pratelli [9] (related to the problems of super-replication and utility maximization in *large financial markets*, i.e., when a sequence of basic assets is considered) to the case of a bond market model. Some proofs will be omitted or outlined when they are only slight modifications of the proofs given by De Donno, Guasoni and Pratelli [9].

We recall that, when $S$ is a $d$-dimensional semimartingale, which models the evolution of the prices of $d$ securities, a predictable $d$-dimensional process $H$ is called a (self-financing) $x$-admissible strategy (according to the terminology of Delbaen and Schachermayer [12]) if it is $S$-integrable and $(H \cdot S)_t \geq -x$. An important consequence of this definition is that, if $S$ is a local martingale, $H \cdot S$ is still a local martingale (by Corollary 3.5 of [2]) and this property is fundamental to study the utility maximization problem.

We have seen in Example 2.1 that, if $\mathbf{S}$ is infinite-dimensional, this property is lost. Conversely, it can be easily proved, by means of the Fatou's lemma, that, if $H^n$ is a sequence of simple integrands such that $H^n \cdot \mathbf{S}$ is bounded from below and $(H^n \cdot \mathbf{S})$ is a Cauchy sequence in $\mathcal{S}(\mathbb{P})$, then the limit process $\mathbf{H} \cdot \mathbf{S}$ is a supermartingale.

Therefore, in the case of a bond market model (and in analogy with De Donno, Guasoni and Pratelli [9]), we give a new definition of admissibility which exploits this fact:

DEFINITION 3.2. Let $x > 0$: a generalized strategy $\mathbf{H}$ is $x$-admissible if there exists an approximating sequence $(H^n)_{n \geq 1}$ of elementary $x$-admissible strategies, such that $(H^n \cdot \overline{\mathbf{P}}) \to (\mathbf{H} \cdot \overline{\mathbf{P}})$ in the semimartingale topology. We denote the set of $x$-admissible generalized strategies by $\mathcal{A}_x$.

We say that $\mathbf{H}$ is *admissible* if it is $x$-admissible for some $x > 0$.

This definition of admissibility allows to use the arguments based on Fatou's lemma, as, for instance, the Ansel and Stricker's theorem (hence, to extend the results by Kramkov and Schachermayer [17]): in particular, if $\mathbf{H}$ is admissible according to the previous definition, then the process $(\mathbf{H} \cdot \overline{\mathbf{P}})_t$ is a supermartingale for every $\mathbb{P}^* \in \mathcal{M}_e$.

In order to apply the convex duality methods to the *utility maximization problem*, the first step is a suitable characterization of super-replicable claims: we have the following result, whose proof is omitted since it is essentially an adaptation of Theorem 3.1 of [9] to the case of a bond market.

THEOREM 3.1. Let $X \in L^0_+$ and $x > 0$. The following conditions are equivalent:

(1) $\sup\limits_{\mathbb{P}^* \in \mathcal{M}_e} E_{\mathbb{P}^*}[X] \leq x$;

(2) There exists $\mathbf{H} \in \mathcal{A}_x$, such that

$$X \leq x + (\mathbf{H} \cdot \overline{\mathbf{P}})_T.$$

Consider a utility function $U: \mathbb{R} \to \mathbb{R}$: we assume that $U(x) = -\infty$ for $x < 0$, that is, negative wealth is not allowed. Furthermore, we assume that $U$ satisfies the Inada (regularity) conditions: $U$ is strictly increasing, strictly concave, continuously differentiable with $U'(0) = \infty$ and $U'(\infty) = 0$.

Let $J$ be the directed set of finite subsets of $I$, ordered by inclusion. For $j = (x_1, \ldots, x_k)$, we denote by $\mathcal{H}^j$ the set of admissible integrands of the form $H = \sum_{x_s \in j} h^s \delta_{x_s}$ and by $u_j(x)$ the maximal expected utility obtained with the bonds $P(\cdot, x_1), \ldots, P(\cdot, x_k)$, namely,

$$u_j(x) = \sup_{H \in \mathcal{H}^j} \mathbb{E}\left[ U\left( x + \int_0^T H_s \, d\overline{\mathbf{P}}_s \right) \right].$$

It is evident that the net of functions $(u_j)_{j \in J}$ is increasing: we pose

$$u_\infty(x) = \lim_{j \in J} u_j(x) = \sup_{j \in J} u_j(x).$$

Moreover, we define the maximal expected utility over all the admissible strategies:

$$u(x) = \sup_{\mathbf{H} \in \mathcal{A}_x} \mathbb{E}\left[ U\left( x + \int_0^T \mathbf{H}_s \, d\overline{\mathbf{P}}_s \right) \right].$$

We have trivially that $u(x) \geq u_\infty(x)$: we shall prove that this is, in fact, an equality. More precisely, for every $j \in J$, we introduce the set

$$\mathcal{C}_j = \{X \in L^0_+ : X \leq 1 + (H \cdot \overline{\mathbf{P}})_T, H \in \mathcal{H}^j\},$$

in such a way that we have

$$(3.1) \qquad u_j(x) = \sup_{X \in \mathcal{C}_j} \mathbb{E}[U(xX)].$$

Let $\mathcal{D}_j$ be the polar of $\mathcal{C}_j$: we recall that the polar of a set $A \subset L^0_+$ is defined by

$$A^\circ = \{f \in L^0_+ : \mathbb{E}[fg] \leq 1 \text{ for all } g \in A\}.$$

It is immediate to verify that the family of sets $\mathcal{D}_j$ is decreasing and that $\mathcal{D} = \bigcap_{j \in J} \mathcal{D}_j$ is the polar of $\mathcal{C} = \bigcup_{j \in J} \mathcal{C}_J$.

As usual, we denote by $V$ the convex conjugate function of $U$, namely,

$$V(y) = \sup_{x > 0}[U(x) - xy]$$

for $y > 0$: the dual problem of (3.1) is then given by

$$(3.2) \qquad v_j(y) = \inf_{Y \in \mathcal{D}_j} \mathbb{E}[V(yY)].$$

The net of functions $(v_j)_{j \in J}$ is increasing, so, in analogy with what we did above with $U$, we can define $v_\infty(y) = \lim_j v_j(y)$ and $v(y) = \inf_{y \in \mathcal{D}} \mathbb{E}[V(yY)]$.

The equality $u(x) = u_\infty(x)$ will be obtained as a consequence of the dual equality $v(y) = v_\infty(y)$ (which is easier to prove). In order to prove the latter result, we have to extend to the case of a net of functions an important result proved by Schachermayer [24] for the case of a sequence of functions (Lemma 3.5). This result is in some sense a *substitute of compactness* in $L^0_+$.

**Lemma 3.2.** *Let $(Z_j)_{j \in J}$ be a net of positive measurable functions, and let $\Gamma_j = \text{Conv}(Z_h | h \geq j)$ be the convex envelope of the set of functions $\{Z_h | h \geq j\}$. It is possible to choose a net of functions $(W_j)_{j \in J}$ and a function $W$ with values in $[0, +\infty]$ in such a way that:*

(1) $W_j \in \Gamma_j$, *for every $j$;*
(2) *the net $(W_j)$ converges to $W$ in probability.*

*Moreover, there exists an increasing sequence $j_1 \leq j_2 \leq \cdots$ such that $W = \lim_n W_{j_n}$ a.s.*

The proof of this lemma (which is a natural extension of the proof given in [12], Lemma A.1.1) can be found in [23]. We can now prove the following result:

LEMMA 3.3.   $v(y) = v_\infty(y)$, for all $y > 0$.

PROOF.   We have evidently $v(y) \geq v_\infty(y)$ for every $y > 0$. Given $j \in J$, let $Y_j \in \mathcal{D}_j$ a function such that $v_j(y) = \mathbb{E}[V(yY_j)]$ (see [17] for the existence of such $Y_j$). Then, we have $v_\infty(y) = \lim_j \mathbb{E}[V(yY_j)]$.

Let us consider $Z_j \in \Gamma_j \subset \mathcal{D}_j$ such that the net $Z_j$ converges to $Z$ in probability: $Z \in \bigcap_j \mathcal{D}_j = \mathcal{D}$. Since $V$ is convex, we have that $\sup_j \mathbb{E}[V(yZ_j)] = v_\infty(y)$.

Let us consider the sequence $(Z_{j_n})_{n \geq 1}$ converging to $Z$ a.s.: since the sequence of negative parts $V^-(yZ_{j_n})$ is uniformly integrable (see Lemma 3.4 in [17] for details), we have

$$v(y) \leq \mathbb{E}[V(yZ)] \leq \liminf_{n \to \infty} \mathbb{E}[V(yZ_{j_n})] = v_\infty(y). \qquad \square$$

It is now easy to adapt the proofs of Lemmas 4.2 and 4.3 in [9] and we have the following result:

PROPOSITION 3.4.   $u_\infty(x) = u(x)$, for all $x > 0$.

With these results, we can immediately give an analogous statement of the main result of Kramkov and Schachermayer [17] (see also Theorem 4.4 of [9]) which shows that the functions $u(\cdot)$ and $v(\cdot)$ are conjugate and, if there exists the optimal solution $\hat{X}(x)$ of the primal problem (3.1), we have the representation $U'(\hat{X}(x)) = \hat{Y}(y)$, where $y = u'(x)$ and $\hat{Y}(y)$ is the optimal solution of the dual problem.

REMARK 3.1.   As in the finite-dimensional case, when a unique martingale measure $\mathbb{P}^*$ does exist, the dual problem takes a very simple form, since it reduces to $v(y) = \mathbb{E}[V(yd\mathbb{P}^*/d\mathbb{P})]$. A special case of this result, together with the use of convex duality methods in the bond market, is illustrated in [8] in the setting of a complete bond market based on a Wiener sheet.

**4. A comparison with other infinite-dimensional approaches.**  In literature one can find several papers devoted to the analysis of the term structure of interest rates by means of infinite-dimensional stochastic integration. We recall, among these, Björk, Di Masi, Kabanov and Runggaldier [4], Filipovic [16], Carmona and Tehranchi [6], Ekeland and Taflin [14] and Aihara and Bagchi [1]. All these authors, except Björk, Di Masi, Kabanov and Runggaldier [4] (whose paper will be discussed in the last part of this section), assume that the bond price process is driven either by an infinite-dimensional Wiener process or by a cylindrical Wiener process. Therefore, a stochastic integral with respect to such processes is needed; expositions of the theory of stochastic integration with respect to an infinite-dimensional (or cylindrical) Wiener process can be found, for instance, in [7] or in [5].

The thesis of Filipovic [16] contains (Chapter 4: *The HJM methodology revisited*) the extension, in an infinite-dimensional setting, of the conditions on the drift in the Heath–Jarrow–Morton model; Carmona and Tehranchi [6] and Aihara and Bagchi [1] consider the question of hedging portfolios for interest rate contingent claims, whereas Ekeland and Taflin [14] study the problem of utility maximization in a bond market. These are essentially the same topics which we analyzed in Section 3, with a different approach.

All these works have in common the fact that the zero coupon bond price process $P(t, \cdot)$ [or the forward rate process $f(t, \cdot)$] is modeled as a stochastic process with values in a Hilbert space $H$ (usually, $H$ is an appropriate weighted Sobolev space), contained in the set of continuous functions on $[0, T^*]$ (or, possibly, $[0, +\infty))$, and with the additional assumption that the "evaluation functional" $\langle \delta_s, g \rangle = g(s)$ belongs to $H'$. This approach is sometimes referred as a "Hilbert space" approach.

We describe in more detail the technique adopted, for instance, by Carmona and Tehranchi [6] in order to compare the "Hilbert space" approach with ours. Carmona and Tehranchi introduce the weighted Sobolev space $F_w^2$, defined as the set of functions $x : \mathbb{R}_+ \to \mathbb{R}$, which are differentiable, such that the derivative $x'$ is absolutely continuous, $x'(\infty) = 0$ and $\int_0^\infty x'(s)^2 w(s) \, ds < \infty$, where $w$ is a given function which represents the "weight." The space $F_w^2$ is a Hilbert space for the norm $\|x\|_{F_w^2} = (\int_0^\infty x'(s)^2 w(s) \, ds < \infty)^{1/2}$. Moreover, they fix a separable Hilbert space $H$ and consider the space of Hilbert–Schmidt operators taking $H$ into $F_w^2$, which is denoted by $\mathcal{L}_{HS}(H, F_w^2)$, and is itself a Hilbert space.

Then, they choose as "state variable" the discounted bond price curve $\overline{\mathbf{P}}_t = \overline{P}(t, \cdot)$ [which in their paper is denoted by $\widetilde{P}_t(\cdot)$]. They assume that, under the risk-neutral measure, the dynamics of $\overline{\mathbf{P}}$ evolves according to the equation $d\overline{\mathbf{P}}_t = \sigma_t \, dW_t$, where $W$ is a cylindrical Wiener process on $H$ and $\sigma_t$ is a process with values in $\mathcal{L}_{HS}(H, F_w^2)$, such that $\int_0^t \|\sigma_s\|_{\mathcal{L}}^2 \, ds < \infty$ a.s. Hence, the process $\overline{\mathbf{P}}$ takes values in the Hilbert space $F_w^2$.

A strategy is defined as a process with values in $(F_v^1)'$, which is the dual space of an appropriate weighted Sobolev space $F_v^1$, chosen in such a way that $F_w^2$ can be continuously embedded in $F_v^1$ (see [6], page 1275 for details). Hence, a strategy takes values in a dense subset of $(F_w^2)'$. In particular, a strategy $\phi$ is defined as a $(F_v^1)'$-valued process such that $\phi(\omega, t)$ belongs to $\overline{\text{span}\{\delta_s; s \geq t\}}^{(F_v^1)'}$ [i.e., the closure of the span$\{\delta_s; s \geq t\}$ in the topology of $(F_v^1)'$] for almost every $(\omega, t)$. This means that each strategy can be approximated with a sequence of stochastic processes, which are finite combinations of Dirac measures, that is, with our definition, simple integrands.

Since the space of continuous local martingales is closed in the semi-martingale topology, it is clear that the strategies considered by Carmona and Tehranchi are, in fact, generalized integrands. Moreover, a generalized integrand which takes values in $(F_v^1)'$ is a strategy according to the definition of Carmona and Tehranchi. The advantage of Carmona and Tehranchi's definition of strategy is that it is particularly easy to give a condition for a strategy to be self-financing (Definition 3.5 in [6]), by means of the duality between Hilbert spaces. Vice versa, we have already pointed out the difficulties arising with our definition of self-financing strategy in the discussion which follows Definition 3.1.

On the other hand, the set of strategies considered by Carmona and Tehranchi is not sufficiently large: even in the case of uniqueness of the martingale measure, the market is only approximately complete (using the terminology introduced by Björk, Di Masi, Kabanov and Runggaldier [4]), that is, the set of hedgeable claims is dense in the set of all sufficiently integrable claims.

On the contrary, our definition of generalized strategy allows to define completeness in the classical sense: more precisely, if the martingale measure is unique, every sufficiently integrable contingent claim can be hedged by means of a generalized integrand.

However, Carmona and Tehranchi [6] show (by using Malliavin's calculus) that if they restrict to a smaller class of claims, a hedging portfolio can be constructed. In particular, they consider *Lipschitz* claims, that is, claims of the form $g(P(T, T_1), \dots, P(T, T_n))$, where $g$ is a Lipschitz function: in this case, an infinite-dimensional version of the Clark–Ocone formula provides an explicit expression for a hedging strategy and conditions are given for the uniqueness of this strategy.

A completely different approach is due to Björk, Di Masi, Kabanov and Runggaldier [4]: they consider the bond prices process as a process with values in a space of continuous functions. With the aim of providing a mathematical background for the theory of bond markets, Björk, Di Masi, Kabanov and Runggaldier [4] suggested the construction of a stochastic integral with respect to processes taking values in a space of continuous functions. We shall prove that their definition is a special case of ours.

Using the same notation of Björk, Di Masi, Kabanov and Runggaldier [4], we call $\mathcal{E}_b$ the set of processes of the form

$$(4.1) \qquad \phi(\omega, t) = \sum_{i \le n} \mathbf{1}_{\Gamma_i \times ]t_i, t_{i+1}]}(\omega, t) m_i,$$

where $m_i \in \mathcal{M}(I)$, $0 \le t_1 < t_2 < \cdots < t_{n+1} \le T$, $\Gamma_i \in \mathcal{F}_{t_i}$.

Björk, Di Masi, Kabanov and Runggaldier [4] consider as integrator a stochastic process $P_t$, which takes values in $\mathcal{C}(I)$ and satisfies the following assumptions (Assumptions 2.1 and 2.2 in [4]):

(i) $P$ is weakly regular in the sense that, for all $\mu \in \mathcal{M}(I)$, the process $\mu(P_t) = \int_I P_t(x) \mu(dx)$ is càdlàg $\mathbb{P}$-a.s.

(ii) $P$ is a controlled process in the following sense: there exists a control pair $(\kappa, p)$, where $\kappa$ is a predictable random measure of the form $\kappa(dt, du) = l_t(du)\,dt$ on $[0, T] \times U$, $\mathcal{B}([0, T]) \otimes \mathcal{U})$, where $(U, \mathcal{U})$ is a Lusin space, while $p$ is a real-valued (measurable) function defined on $(\Omega \times [0, T] \times U \times \mathcal{M}(I), \mathcal{P}r \times \mathcal{U} \times \sigma(\mathcal{M}(I)))$ ($\mathcal{P}r$ denotes the predictable $\sigma$-field), with the following properties:

(a) $K_t = 1 + \kappa([0, t] \times U) < \infty$ for all $t$;

(b) $p(\omega, t, u, \cdot)$ is weakly continuous, it is a seminorm on $\mathcal{M}(I)$ such that $p(\omega, t, u, \mu) \leq \|\mu\|_V$ [where $\| \cdot \|_V$ denotes the total variation norm on $\mathcal{M}(I)$];

(c) given a process $\phi \in \mathcal{E}_b$ of the form (4.1) and define the integral $(\phi \cdot P)$ in the natural way as

$$(4.2) \qquad (\phi \cdot P)_t = \sum_{i \leq n} (m_i(P_{t_{i+1} \wedge t}) - m_i(P_{t_i \wedge t})) \mathbf{1}_{\Gamma_i},$$

the following inequality holds for any stopping time $\tau \leq T$:

$$(4.3) \qquad \mathbb{E}\Big[\sup_{t \leq \tau} |(\phi \cdot P)_t|^2\Big] \leq C\mathbb{E}\Big[K_\tau \int_0^\tau \int_U p^2(s, u, \phi_s)\kappa(ds, du)\Big],$$

where $C$ is a constant which does not depend on $\phi$.

With these assumptions, the process $P$ satisfies the assumptions on the integrator process that we made in Section 2. Indeed, if we set $S^x = P_t(x) = (\delta_x \cdot P)_t$, we have that, by condition (i) above, $S^x$ is a càdlàg process. Furthermore, it can be proved, as a particular case of Theorem 2.4(a) in [4], that $S^x$ is a semimartingale. So, we have defined a family of semimartingales $\mathbf{S} = (S^x)_{x \in I}$.

To show that $\mathbf{S}$ fulfills Assumption 2.1, that is, the mapping $x \mapsto S^x$ is continuous, we can follow the proof of Theorem 2.4(a) in [4]. We can assume, by localization, that $\mathbb{E}[K_T^2] < \infty$. Let $H$ be a real predictable elementary integrand, uniformly bounded by 1. Then, setting $\phi = H(\delta_x - \delta_y)$, it follows from (4.3) that

$$\mathbb{E}\Big[\sup_{t \leq T} |(H \cdot S^x)_t - (H \cdot S^y)_t|^2\Big] \leq C\mathbb{E}\Big[K_T \int_0^T \int_U p^2(s, u, \delta_x - \delta_y)\kappa(ds, du)\Big].$$

In particular, the following inequality holds:

$$\|S^x - S^y\|_{\mathcal{S}(\mathbb{P})} \leq C\mathbb{E}\Big[K_T \int_0^T \int_U p^2(s, u, \delta_x - \delta_y)\kappa(ds, du)\Big].$$

Because of condition (b), we have that $p(s, u, \delta_x - \delta_y) \leq \|\delta_x - \delta_y\|_V = 2$; since $\lim_{y \to x} p^2(s, u, \delta_x - \delta_y) = 0$, it follows that $S^y$ converges to $S^x$ in $\mathcal{S}(\mathbb{P})$ as $y$ tends to $x$.

It remains to show that the integrands according to Björk, Di Masi, Kabanov and Runggaldier [4] are generalized integrands. Let $\tau$ be a fixed bounded stopping time. Björk, Di Masi, Kabanov and Runggaldier denote by $L^2_\tau$ the set of all weakly predictable measure-valued processes $\phi$, such that

$$q^2_\tau(\phi) = \mathbb{E}\left[K_\tau \int_0^\tau \int_U p^2(s, u, \phi_s)\kappa(ds, du)\right] < \infty$$

and by $L^2_{\mathrm{loc}}(\mathbb{P})$ the set of all weakly predictable measure-valued processes $\phi$, such that, for all $t$,

$$\int_0^t \int_U p^2(s, u, \phi_s)\kappa(ds, du) < \infty \qquad \text{a.s.}$$

They show that the integral defined on $\mathcal{E}_b$ by (4.2) can be extended to $L^2_\tau$ continuously with respect to the seminorm $\Pi_\tau(Y) = \mathbb{E}[\sup_{t \leq \tau} |Y_t|^2]^{1/2}$ and then, by localization, to $L^2_{\mathrm{loc}}(\mathbb{P})$. Furthermore, they prove that the integral $(\phi \cdot P)$ does not depend on the particular choice of a control pair $(\kappa, p)$ and that it is a semimartingale.

In order to show that every process in $L^2_{\mathrm{loc}}(\mathbb{P})$ is a generalized integrand in the sense of Definition 2.6, it is sufficient to show that every $\phi \in L^2_\tau$ is a generalized integrand for any bounded stopping time $\tau$. Björk, Di Masi, Kabanov and Runggaldier [4] prove that such a process is limit of a sequence $\phi^n \in \mathcal{E}_b$, such that the sequence of integrals $\phi^n \cdot P$ converges in the topology given by $q_\tau$, hence, in $\mathcal{S}(\mathbb{P})$. So, in fact, we only need to prove that, for every process $\phi \in \mathcal{E}_b$, there exists a sequence of simple integrands $H^n$ such that $H^n \cdot \mathbf{S}$ converges to $\phi \cdot \mathbf{S}$ in $\mathcal{S}(\mathbb{P})$. In particular, it is sufficient to show that this holds for a process $\phi$ of the form $\phi = \mathbf{1}_A m$, where $A$ is a predictable set and $m$ is an element of $\mathcal{M}(I)$. It is a known fact that $m$ is the limit (in the weak topology) of a sequence $m_n$ of linear combination of Dirac measures, such that $\|m_n\|_V \leq \|m\|_V$. Then, $\phi_n = \mathbf{1}_A m_n$ is a sequence which satisfies our requirements.

## REFERENCES

[1] Aihara, S. I. and Bagchi, A. (2005). Stochastic hyperbolic dynamics for infinite-dimensional forward rates and option pricing. *Math. Finance* **15** 27–47. MR2116795

[2] Ansel, J. P. and Stricker, C. (1994). Couverture des actifs contingents et prix maximum. *Ann. Inst. H. Poincaré Probab. Statist.* **30** 303–315. MR1277002

[3] Björk, T. (1998). *Arbitrage Theory in Continuous Time.* Oxford Univ. Press.

[4] Björk, T., Di Masi, G., Kabanov, Y. and Runggaldier, W. (1997). Towards a general theory of bond markets. *Finance Stoch.* **1** 141–174.

[5] Carmona, R. (2005). Interest rate models: From parametric statistics to infinite-dimensional stochastic analysis. *SIAM.* To appear.

[6] CARMONA, R. and TEHRANCHI, M. (2004). A charachterization of hedging port-folios for interest rate contingent claims. *Ann. Appl. Probab.* **14** 1267–1294. MR2071423

[7] DA PRATO, G. and ZABCZYK, J. (1992). *Stochastic Equations in Infinite Dimensions.* Cambridge Univ. Press. MR1207136

[8] DE DONNO, M. (2004). The term structure of interest rates as a random field: A stochastic integration approach. In *Stochastic Processes and Applications to Mathematical Finance* 27–51. World Scientific Publishing, River Edge, NJ.

[9] DE DONNO, M., GUASONI, P. and PRATELLI, M. (2005). Super-replication and utility maximization in large financial markets. *Stoch. Process. Appl.* To appear.

[10] DE DONNO, M. and PRATELLI, M. (2004). On the use of measure-valued strategies in bond markets. *Finance Stoch.* **8** 87–109. MR2022980

[11] DE DONNO, M. and PRATELLI, M. (2005). Stochastic integration with respect to a sequence of semimartingales. *Séminaire de Probabilités XXXIX.* To appear.

[12] DELBAEN, F. and SCHACHERMAYER, W. (1994). A general version of the fundamental theorem of asset pricing. *Math. Ann.* **300** 463–520. MR1304434

[13] DUNFORD, N. and SCHWARTZ, J. T. (1988). *Linear Operators*, **1**. Wiley, New York. MR1009162

[14] EKELAND, I. and TAFLIN, E. (2005). A theory of bond portfolios. *Ann. Appl. Probab.* **15**. MR2134104

[15] EMERY, M. (1979). Une topologie sur l'espace des semi-martingales. *Séminaire de Probabilités XIII. Lecture Notes in Math.* **721** 260–280. Springer, Berlin. MR544800

[16] FILIPOVIC, D. (2001). *Consistency Problems for Heath–Jarrow–Morton Interest Rate Models. Lecture Notes in Math.* **1760**. Springer, Berlin. MR1828523

[17] KRAMKOV, D. and SCHACHERMAYER, W. (1999). The asymptotic elasticity of utility functions and optimal investment in incomplete markets. *Ann. Appl. Probab.* **9** 904–950. MR1722287

[18] MÉMIN, J. (1980). Espace de semi-martingales et changement de probabilité. *Z. Wahrsch. Verw. Gebiete* **52** 9–39. MR568256

[19] MÉTIVIER, M. (1982). *Semimartingales (a Course on Stochastic Processes).* de Gruyter, Berlin. MR688144

[20] MIKULEVICIUS, R. and ROZOVSKII, B. L. (1998). Normalized stochastic integrals in topological vector spaces. *Séminaire de Probabilités XXXII. Lecture Notes in Math.* **1686** 137–165. Springer, Berlin. MR1655149

[21] MIKULEVICIUS, R. and ROZOVSKII, B. L. (1999). Martingales problems for stochastic PDE's. In *Stochastic Partial Differential Equations*: *Six Perspectives* (R. Carmona, B. Rozovski, eds.) 243–325. Amer. Math. Soc., Providence, RI. MR1661767

[22] PHAM, H. (2003). A predictable decomposition in an infinite assets model with jumps. Application to hedging and optimal investment. *Stochastics Stochastics Rep.* **75** 343–368. MR2017783

[23] PRATELLI, M. (2005). A minimax theorem without compactness hypothesis. *Mediterr. J. Math.* **2** 103–112. MR2135081

[24] SCHACHERMAYER, W. (1992). A Hilbert space proof of the fundamental theorem of asset pricing in finite discrete time. *Insurance Math. Econom.* **11** 249–257. MR1211972

Dipartimento di Matematica
Università di Pisa
Largo B. Pontecorvo 5
56127 Pisa
Italy
e-mail: mdedonno@dm.unipi.it
        pratelli@dm.unipi.it